\documentclass[10pt]{amsart}
\usepackage{amssymb,latexsym,amsmath}
\input{xypic}
\xyoption{all}

\begin{document}

\begin{abstract} 
This is an announcement of results on rank-one Eisenstein cohomology of ${\rm GL}_N$, with $N \geq 3$ an odd integer, and algebraicity theorems for ratios of successive critical values of certain Rankin--Selberg $L$-functions for ${\rm GL}_n \times {\rm GL}_{n'}$ when $n$ is even and $n'$ is odd.  \\
R\'esum\'e: Cette note est une annonce de r\'esultats sur la cohomologie d'Eisenstein de rang un de $ {\rm GL} _ N$,  
avec $N \geq 3$ un entier impair, et de 
th\'eor\`emes  d'alg\'ebricit\'e pour les rapports de valeurs critiques successives de certaines fonctions $L$ de Rankin--Selberg pour 
${\rm GL} _n \times {\rm GL} _ {n '}$ lorsque $n$ est pair et $n'$ est impair.
\end{abstract}

\title[Ratios of critical values]{Eisenstein Cohomology and ratios of critical 
values of Rankin--Selberg $L$-functions} 

\author{\bf G\"unter Harder \ \ \and \ \ A. Raghuram}

\date{\today}
\thanks{Work of A.R. is partially supported by NSF grant DMS-0856113 and an Alexander von Humboldt Research Fellowship}

\maketitle


\def\g{\mathfrak{g}}
\def\k{\mathfrak{k}}
\def\z{\mathfrak{z}}
\def\l{\mathfrak{l}}
\def\gl{\mathfrak{gl}}
\def\w{{\bf w}}
\def\Ext{{\rm Ext}}
\def\Aut{{\rm Aut}}
\def\Hom{{\rm Hom}}
\def\Ind{{\rm Ind}}
\def\Id{{\rm Id}}
\def\GL{{\rm GL}}
\def\SL{{\rm SL}}
\def\SO{{\rm SO}}
\def\O{{\rm O}}
\def\Gal{{\rm Gal}}
\def\End{{\rm End}}
\def\R{\mathbb{R}}
\def\C{\mathbb{C}}
\def\Z{\mathbb{Z}}
\def\Q{\mathbb{Q}}
\def\A{\mathbb{A}}
\def\M{\bf M}
\def\Hp{H_{\pi_f}}
\def\Hpp{H_{\pi_p}}
\def\HCl{H_{\pi_\infty^\lambda}}

\def\CH{\mathcal{H}}
\def\CM{\mathcal{M}}
\def\CA{\mathcal{A}}
\def\CO{\mathcal{O}}
\def\CI{\mathcal{I}}
 
\def\bfpi{\mathbf{\Pi}}
\def\bfdelta{\mathbf{\Delta}}

\def\Cat{\mathcal{C}}
\def\E{\mathcal{E}}

\def\atimes{\stackrel{\rm a}{\times}}

\def\to{\rightarrow}
\def\To{\longrightarrow}
\def \iso{ \buildrel \sim \over\longrightarrow}
\def\pkt{\bullet}
\def\into{\hookrightarrow}

\def\1{1\!\!1}
\def\dim{{\rm dim}}

\def\th{^{\rm th}}
\def\isom{\approx}
 \def\rn{^\circ r}



\numberwithin{equation}{section}
\newtheorem{thm}[equation]{Theorem}
\newtheorem{cor}[equation]{Corollary}
\newtheorem{lemma}[equation]{Lemma}
\newtheorem{prop}[equation]{Proposition}
\newtheorem{con}[equation]{Conjecture}
\newtheorem{ass}[equation]{Assumption}
\newtheorem{defn}[equation]{Definition}
\newtheorem{rem}[equation]{Remark}
\newtheorem{exer}[equation]{Exercise}
\newtheorem{exam}[equation]{Example}

\section{The general situation}
Let $G/\Q$ be a connected split reductive algebraic group over $\Q$ whose derived group $G^{(1)}/\Q$ is simply connected. 
Let $Z/\Q$ be the center of $G$ and let $S$ be the maximal $\Q$-split torus in $Z$. 
Let $C_{\infty}$ be a maximal compact subgroup of $G(\R)$ and let $K_\infty = C_\infty S(\R)^\circ$. 
The connected component of the identity of $K_{\infty}$ is denoted $K_{\infty}^\circ$ and 
$K_{\infty}/K_{\infty}^\circ = \pi_0(K_{\infty}) \iso \pi_0(G(\R)).$
Let $K_f = \prod_p K_p \subset G(\A_f)$ be an open compact subgroup; here $\A$ is the ad\`ele ring of $\Q$ and 
$\A_f$ is the ring of finite ad\`eles. The locally symmetric space of $G$ with level structure $K_f$ is defined as 
$$ 
S^G_{K_f} := G(\Q)\backslash G(\A)/K_{\infty}^{\circ}K_f.
$$
(For the following see Harder~\cite[Chapter 3, Sections 2, 2.1, 2.2] {harder-book} for details.) 
For a dominant integral weight $\lambda$, let $E_{\lambda}$ be an absolutely irreducible finite-dimensional representation of $G/\Q$ 
with highest weight $\lambda$, 
and let $\E_{\lambda}$ denote the associated sheaf on $S^G_{K_f}$. 
We have an action of the Hecke-algebra $ \CH=\CH_{K_f}^G = \otimes_p' \CH_p$ on the cohomology groups $H^{\bullet}(S^G_{K_f}, \E_{\lambda}).$ 

We always fix a level, but sometimes drop it in the notation.
For any finite extension $F/\Q,$ let $E_{\lambda, F} = E_{\lambda} \otimes_\Q F$, then $\E_{\lambda, F}$ is the corresponding sheaf on $S^G_{K_f}$.

Let $\bar{S}^G_{K_f}$ be the Borel--Serre compactification of $S^G_{K_f} $, i.e., 
$\bar{S}^G_{K_f} = S^G_{K_f} \cup \partial \bar{S}^G_{K_f}$, where the boundary is stratified as 
$\partial \bar{S}^G_{K_f}  = \cup_P \partial_PS^G_{K_f}$ with $P$ running through the conjugacy classes of proper parabolic subgroups 
defined over $\Q$. The sheaf $\E_{\lambda, F}$ on $S^G_{K_f}$ naturally extends, using the definition of the Borel-Serre compactification, to a sheaf on 
$\bar{S}^G_{K_f}$ which we also denote by $\E_{\lambda, F}$. 
Restriction from $\bar{S}^G_{K_f}$ to $S^G_{K_f}$ in cohomology induces an isomorphism $H^i(\bar{S}^G,\E_{\lambda}) \iso  H^i(S^G,\E_{\lambda}).$

Our basic object of interest is the following long exact sequence of $\pi_0(K_{\infty}) \times \CH$-modules 
$$
\cdots  \longrightarrow H^i_c(S^G, \E_{\lambda}) \stackrel{\iota^*}{\longrightarrow}   H^i(\bar{S}^G,\E_{\lambda}) 
\stackrel{r^*}{\longrightarrow } H^i(\partial \bar{S}^G,\E_{\lambda}) 
\longrightarrow H^{i+1}_c(S^G, \E_{\lambda}) \longrightarrow \cdots
$$
The image of cohomology with compact supports inside the full cohomology is called {\it inner} or {\it interior} cohomology and is denoted 
$H^{\bullet}_{\, !} := {\rm Image}(\iota^*) = {\rm Im}(H^{\bullet}_c \to H^{\bullet})$. The theory of Eisenstein cohomology is designed to describe the image of 
the restriction map $r^*$. Our goal is to study the arithmetic information contained in the above exact sequence.  

The inner cohomology  is a semi-simple module for the Hecke-algebra. (See Harder \cite[Chap. 3, 3.3.5]{harder-book}.) After a suitable finite extension $F/\Q,$   where $\Q\subset F\subset \bar\Q \subset \C,$ we have an isotypical decomposition 
$$  
H^i_{\,!}(S^G_{K_f}, \E_{\lambda,F}) =\bigoplus_{\pi_f\in {\rm Coh}(G,K_f,\lambda)} H^i_{\,!}(S^G_{K_f}, \E_{\lambda,F})(\pi_f)
$$
where $\pi_f$ is an isomorphism type of an absolutely  irreducible $\CH$-module, i.e., an $F$-vector space $H_{\pi_f}$
with an absolutely irreducible action of $\CH$.  The local factors  $\CH_p$ are commutative outside a finite set $V=V_{K_f}$ of primes
 and the factors $\CH_p$ and $\CH_q$, for $p\neq q$, commute with each other.  Hence for $p\not\in V$ the commutative algebra $\CH_p$ acts
on $H_{\pi_f} $  by a homomorphism $\pi_p:\CH_p \to F.$ Let $H_{\pi_p}$ be the one dimensional vector space $F$ with basis $1\in F$ 
with the action $\pi_p$ on it. Then  $H_{\pi_f}  =\otimes_{p\in V}H_{\pi_p}\otimes^\prime_{p\not\in V}H_{\pi_p}=\otimes^\prime_p H_{\pi_p}.$
The set of isomorphism classes which occur in  the  above decomposition is called the `spectrum' ${\rm Coh}(G, K_f, \lambda)$.
 If we restrict the elements of the Galois group  $\Gal(\bar{\Q}/\Q)$ to $F$  we get the conjugate embeddings of $F$ into $\bar\Q;$ 
we introduce ${\CI}(F) =\{\iota : F\to \C\}=\{\iota : F\to \bar\Q\}.$   For $\iota \in \CI(F)$ define $\iota\circ \pi_f$ as $\Hp\otimes_{F,\iota}\C$.  
We define the rationality field of $\pi_f$ as 
$\Q(\pi_f)=\{x\in F  \, \vert \, \iota(x)=\iota^\prime(x)  \text{ if } \iota\circ\pi_f=\iota^\prime\circ \pi_f\}.$

\section{The case of $\GL_n$ and the definition of relative periods when $n$ is even}

 Let $T /\Q $ be a maximal $\Q$-split torus in $G$, let $T^{(1)}=T\cap G^{(1)}.$ Let $X^*(T)$ be its group of characters then restriction
 of characters gives an inclusion  $X^*(T)\subset  X^*(T^{(1)}) \oplus X^*(Z)$ and after tensoring by $\Q$ this becomes an isomorphism.
Any $\lambda\in X^*(T) $ can be written as $\lambda^{(1)}+\delta, \lambda^{(1)}\in  X^*(T^{(1)})\otimes \Q=:X_\Q^*(T^{(1)}),\delta\in X^*_\Q(Z).$ 

Consider the case  $G = \GL_n/\Q.$ Take a regular essentially self-dual dominant integral highest weight $\lambda$. Let $\rho\in X^*_\Q(T^{(1)}) $ be half the sum of positive roots, and write 
$ \lambda+\rho  = a_1\gamma_1+\dots +a_{n-1} \gamma_{n-1}+d\cdot {\rm det},$ which is an equation in $X_\Q^*(T)
$; 
the  $\gamma_i\in  X_\Q^*(T) $ restrict to the fundamental weights in $X^*(T^{(1)})$ 
and are trivial on the center $Z$. 
Regular, dominant and integral mean that $a_i\geq 2$ are integers, and essentially self-dual means $a_i=a_{n-i}$. Further, 
for such a weight $\lambda$ we have $2d \in \Z$ and it satisfies the parity condition: 
\begin{equation}
\label{eqn:parity}
2d \equiv \w + n - 1 \pmod{2}
\end{equation}
where $\w = \w(\lambda) := \sum_i a_i$ is the `motivic weight'; see below.

Given such a $\lambda$, there is a unique essentially unitary Harish-Chandra module $\HCl$  such that the relative Lie algebra cohomology group 
$H^\bullet(\g,K_\infty^\circ, \HCl \otimes E_\lambda)\neq  0.$ Let $L^2_d(G(\Q)\backslash G(\A)/K_f, \omega_{E_{\lambda}}^{-1})$ denote the 
discrete spectrum for $G(\A)$ in the space of $L^2$-automorphic forms with level structure $K_f$ on which $Z(\R)^{\circ}$ acts via the inverse of 
the central character of $E_{\lambda}$. For $\pi_f \in {\rm Coh}(G,K_f,\lambda)$ and $\iota \in \CI(F)$   we consider
$$
W(\pi_\infty^\lambda\otimes\iota \circ \pi_f) = 
 \Hom_{(\g, K_\infty^\circ)\times \CH^G_{K_f}} \left(\HCl\otimes (\Hp\otimes_{F,\iota}\C), \, L^2_d(G(\Q)\backslash G(\A)/K_f, \omega_{E_{\lambda}}^{-1}) \right)
$$  
which is one-dimensional due to multiplicity-one for the discrete spectrum of $\GL_n$; the image is in fact in the cuspidal spectrum by regularity of $\lambda$. 
(See, for example, Schwermer \cite[Corollary 2.3]{schwermer}.) 
We choose a generator $\Phi$ for $W(\pi_\infty^{\lambda} \times \iota \circ \pi_f ).$

The summand $H^{\bullet}_!(S^G_{K_f}, \E_{\lambda , F}) (\pi_f)$ can be decomposed for the action of $\pi_0(G(\R)) = \Z/2\Z$ as
$$  
H^{\bullet}_!(S^G_{K_f}, \E_{\lambda,F}) (\pi_f)=\bigoplus_{\epsilon :\pi_0(G(\R))\to \Z/2\Z}
H^{\bullet}_!(S^G_{K_f}, \E_{\lambda,F}) (\pi_f)(\epsilon).
$$
The action of  $\pi_0(G(\R)) =  \pi_0(K_{\infty}) = K_{\infty}/K_{\infty}^\circ$ is via its action on 
$H^\bullet(\g,K_\infty^\circ, H_{\pi^\lambda_\infty}\otimes E_{\lambda}).$ (See, for example, Borel-Wallach \cite[I.5]{borel-wallach}.)  
Therefore, we get 
$$ 
\bigoplus_{\epsilon}W(\pi^\lambda_\infty\otimes \iota \circ \pi_f) \otimes H^\bullet(\g,K_\infty^\circ, H_{\pi^\lambda_\infty}\otimes E_{\lambda})(\epsilon)
 \otimes H_{\pi_f} \otimes_{F,\iota}\C \to  
 \bigoplus_{\epsilon}H^{\bullet}_{!}(S^G_{K_f}, \E_{\lambda, F}) (\pi_f)\otimes_{F,\iota} \C(\epsilon).
 $$
Let $b_n = n^2/4$ if $n$ is even, and $(n^2-1)/4$ if $n$ is odd. Since $\pi$ is cuspidal,  it is well-known (see, for example, Clozel \cite{clozel}) that 
$\pi^\lambda_\infty$ is irreducibly induced from essentially discrete series representations and that
 $$
 H^{b_n}(\g,K_\infty^\circ, H_{\pi^\lambda_\infty}\otimes E_\lambda)= 
 \begin{cases}  
 H^{b_n}(\g,K_\infty^\circ, H_{\pi^\lambda_\infty}\otimes E_\lambda)_+ \oplus 
 H^{b_n}(\g,K_\infty^\circ, H_{\pi^\lambda_\infty}\otimes E_\lambda)_- & \mbox{if $n$  is  even;} \cr
H^{b_n}(\g,K_\infty^\circ, H_{\pi^\lambda_\infty}\otimes E_\lambda)_{\epsilon}& \mbox{if $n$ is odd,}
\end{cases}
$$		
where each piece on the right hand side is one-dimensional, and $\epsilon$ is a canonical sign (see \cite[Section 3.3]{raghuram-shahidi-imrn}). 

Now let $n$ be even. We  will define certain periods that we call {\it relative periods}.  We define consistent choices of generators
     $$ \omega_+\in \Hom_{K_\infty^\circ}(\Lambda^{b_n}(\g/\k), H_{\pi^\lambda_\infty}\otimes E_\lambda)_+, \ \ 
      \omega_-\in \Hom_{K_\infty^\circ}(\Lambda^{b_n}(\g/\k), H_{\pi^\lambda_\infty}\otimes E_\lambda)_-, 
      $$
from which we get isomorphisms 
 $$(\Phi\otimes \omega_{\pm}) :   
  \iota \circ \pi_f   \to H^{b_n}_!(S^G_{K_f}, \E_{\lambda,F}) (\pi_f)_{\pm} \otimes_{F,\iota} \C.  
$$  
Composing the inverse of one with the other gives a canonical transcendental isomorphism 
\begin{equation}
\label{eqn:trans}
T^{\rm trans}(\pi_f,\iota) = (\Phi\otimes \omega_-) \circ (\Phi\otimes \omega_+)^{-1}
: H^{b_n}_!(S^G_{K_f}, \E_{\lambda, F}) (\pi_f)_+
\otimes_{F,\iota} \C\to H^{b_n}_!(S^G_{K_f}, \E_{\lambda, F}) (\pi_f)_-
\otimes_{F,\iota} \C. 
\end{equation}
This isomorphism does not depend on the choice of $\Phi$ or the pair $(\omega_+,\omega_-)$ because these are unique up to scalars
which cancel out.  On the other hand we have an arithmetic isomorphism of $\CH^G_{K_f}$-modules 
\begin{equation}
\label{eqn:arith}
T^{\rm arith}(\pi_f):  H^{b_n}_!(S^G_{K_f}, \E_{\lambda, F}) (\pi_f)_+\to
 H^{b_n}_!(S^G_{K_f}, \E_{\lambda, F}) (\pi_f)_- 
\end{equation}
which is unique up to an element in $\Q(\pi_f)^\times.$ Comparing (\ref{eqn:trans}) with (\ref{eqn:arith}) we get the following definition. 

\begin{defn}
There is an array of complex numbers $ \Omega(\pi_f)  =  (\dots, \Omega(\pi_f,\iota), \dots )_{\iota\in\CI(F)}$ defined by 
$$ \Omega(\pi_f,\iota)T^{\rm trans}(\pi_f,\iota)=T^{\rm arith}(\pi_f)\otimes_{F,\iota}\C .$$
 Changing $T^{\rm arith}(\pi_f)$ by an element $a\in \Q(\pi_f)^\times $  changes the array into $ (\dots, \Omega(\pi_f,\iota)\iota(a), \dots )_{\iota:F\to \C}.$
\end{defn}

If we pass from  $\lambda$ to  $\lambda- l \cdot {\rm det}$ for an integer $l$,  then we have a canonical isomorphism
$$  
H^{\bullet}_!(S^G_{K_f}, \E_{\lambda, F}) (\pi_f) \to   
H^{\bullet}_!(S^G_{K_f}, \E_{\lambda - l \cdot {\rm det}, F}) (\pi_f\otimes \vert \ \vert^l)
$$
under which the $\pm $ components are switched by $(-1)^l.$ We get the following period relation: 
\begin{equation}
\label{eqn:period-relations}
\Omega(\pi_f,\iota) =  \Omega(\pi_f\otimes \vert \ \vert^l,\iota)^{(-1)^l}.
\end{equation} 

\begin{rem}{ \rm
Since cuspidal automorphic representations of $\GL_n$ are globally generic we can also define periods by comparing rational structures on 
Whittaker models and cohomological realizations. The periods were denoted $p^{\pm}(\pi_f)$ in Raghuram--Shahidi \cite{raghuram-shahidi-imrn} and 
they appear in algebraicity results for the central critical value of Rankin--Selberg $L$-functions for $\GL_n \times \GL_{n-1}$; 
see Raghuram~\cite[Theorem 1.1]{raghuram-imrn}. The periods $p^{\pm}(\pi_f)$ depend on a choice of a nontrivial character of $\Q \backslash \A$ which is implicit in any discussion concerning Whittaker models. However, one may check that if we change this character then the period changes only by an element of  $\Q(\pi_f)^\times .$  Further, it is an easy exercise to see that $\Omega(\pi_f) = p^+(\pi_f)/p^-(\pi_f)$ up to elements in 
$\Q(\pi_f)^\times.$ On the other hand, the definition of the relative periods $\Omega(\pi_f)$ does not require Whittaker models suggesting that it is far more intrinsic to the representation viewed as a Hecke-summand of global cohomology.}
\end{rem}

\section{The case $G=\GL_n\times \GL_{n^\prime}$ with $n$ even and $n'$ odd}

Let $\sigma_f\in {\rm Coh}(\GL_n, \lambda)$ and $\sigma'_f \in {\rm Coh}(\GL_{n'}, \lambda').$ 
The level structures will be suppressed from our notation from now on. As before, the weights are written as 
$ \lambda+\rho  = a_1\gamma_1+\dots +a_{n-1} \gamma_{n-1} +d \cdot {\rm det},$ and similarly 
$\lambda^\prime+\rho^\prime =  a^\prime_1\gamma^\prime_1+\dots +a^\prime_{n^\prime-1} \gamma^\prime_{n^\prime-1}+d^\prime \cdot {\rm det}^\prime,$
where $a_i=a_{n-i},$ $a^\prime_{i}=a^\prime_{n^\prime-i}$, and again we assume regularity for both the weights.   Let $G=\GL_n\times \GL_{n^\prime}$, $\mu =  \lambda+ \lambda^\prime$ and $\pi_f = \sigma_f \times \sigma'_f$. 
By the K\"unneth formula we get 
$$    
H^{\bullet}_!(S^G, \E_{\mu, F}) (\pi_f) =  
H^{\bullet}_!(S^{\GL_n}, \E_{\lambda^\prime, F}) (\sigma_f) \otimes 
H^{\bullet}_!(S^{\GL_{n'}}, \E_{\lambda^\prime, F}) (\sigma^\prime_f).
$$

Using Grothendieck's conjectural theory of motives, one supposes 
that there are motives ${\bf M}_{\rm eff}$ (resp., ${\bf M}'_{\rm eff}$) that are conjecturally attached to $\sigma_f$ (resp., $\sigma'_f$). (See, for example,  \cite{motives}.) We call a pair of integers $(p,q)$ a Hodge-pair for a motive ${\bf M}$ if the Hodge number $h^{p,q}({\bf M}) \neq 0$. 
The Hodge-pairs of  the motives ${\bf M}_{\rm eff}$ (resp., ${\bf M}'_{\rm eff}$) are expected to be 
$\{ (\w,0), (\w-a_1,a_1)  , \dots,   (0,\w)\}$ (resp., $\{(\w^\prime ,0), (\w^\prime -a^\prime_1,a^\prime_1), \dots,   (0, \w^\prime)\}$) 
where $\w = \sum_{i=1}^{n-1} a_i$ (resp., $\w' = \sum_{i'=1}^{n'-1} a_{i'}'$) are the motivic weights. The motives ${\bf M}_{\rm eff}$ (resp., ${\bf M}'_{\rm eff}$) are suitable Tate-twists of the motives expected to be attached to $\sigma_f$ (resp., $\sigma'_f$) as in Clozel \cite[Conjecture 4.5]{clozel}. The assertion about Hodge pairs may be verified by working with the representations at infinity and their associated local $L$-factors which determine the $\Gamma$-factors  at infinity. The set of Hodge-pairs for ${\bf M}_{\rm eff} \otimes {\bf M}_{\rm eff}'$ are all the pairs of the form
$(\w-a_1\dots -a_s+\w^\prime-a^\prime_1 \dots-a^\prime_{s^\prime}, a_1+\dots+a_s+a_1'+\dots+a_{s'}').$
 
The motivic $L$-function $L( {\bf M}_{\rm eff} \otimes {\bf M}'_{\rm eff}, \iota, s)$ is defined as in Deligne~\cite[(1.2.2)]{deligne}. Intimately related to it is a `cohomological' 
$L$-function $L^{\rm coh}(\sigma_f \times \sigma^\prime_f  ,\iota,s)$ which is defined as an Euler product, where each Euler factor is expressed in terms of eigenvalues of certain normalized Hecke-operators acting on integral cohomology groups. Assume that the middle Hodge number of 
${\bf M}_{\rm eff} \otimes {\bf M}_{\rm eff}'$ 
is zero, i.e., $h^{(\w+\w^\prime)/2, (\w+\w^\prime)/2}=0.$ Put 
$p(\mu)  := \min \{p\; \vert \; \w+\w^\prime\geq  p  >(\w+\w^\prime)/2, \ h^{p,\w+\w^\prime-p} \neq 0  \}.$
Let $\sigma^{\prime{\sf v}}$ denote the contragredient of $\sigma^{\prime}.$
The critical points of $L^{\rm coh}(\sigma_f \times \sigma_f^{\prime{\sf v}}  ,\iota,s)$ are the integers 
\begin{equation}
\label{eqn:critical}
\{ p(\mu), \, p(\mu)-1, \,  \dots \, , \, \w+\w'+1-p(\mu)\}.
\end{equation}
Note that this decreasing list of integers is centered around $(\w+\w^\prime+1)/2$ which is the center of symmetry of the cohomological $L$-function. The total number of critical integers is $2p(\mu)-(\w+\w')$.
The cohomological $L$-function is up to a shift in the $s$-variable the usual automorphic Rankin--Selberg $L$-function 
$L(\sigma_f \times \sigma_f^{\prime{\sf v}}, \iota,s) := L\left( ( \iota \circ \sigma_f ) \times (\iota \circ \sigma_f^{\prime{\sf v}}),s \right)$ for which the functional equation is between $s$ and $1-s$. More precisely, we have 
\begin{equation}
\label{eqn:cohomological-automorphic}
L^{\rm coh}\left(\sigma_f \times \sigma_f^{\prime{\sf v}} ,\iota,s \right) = 
L\left(\sigma_f \times \sigma_f^{\prime{\sf v}}  ,\iota, s - \frac{(\w+\w')}{2} + a(\mu) \right) 
\end{equation}
where $a(\mu) = d-d'$. The parity condition (\ref{eqn:parity}) when applied to both the weights $\lambda$ and $\lambda'$ implies that the shift in 
$- \frac{(\w+\w')}{2} + a(\mu)$ in the $s$-variable is always a half-integer. Observe that the cohomological $L$-function is invariant under changing 
$\sigma$ to $\sigma \otimes \vert \ \vert^l$ or $\sigma'$ to $\sigma' \otimes \vert \ \vert^{l'}.$

A celebrated conjecture of Deligne predicts the existence of two periods $\Omega_{\pm}({\bf M}_{\rm eff} \otimes {\bf M}_{\rm eff}')$ obtained from the Betti and de~Rham realizations of this motive that capture, up to prescribable powers of $(2\pi i)$, the possibly transcendental parts of the critical values of 
$L( {\bf M}_{\rm eff} \otimes {\bf M}_{\rm eff}', \iota, s)$. See \cite[Conjecture 2.7, (3.1.2) and (5.1.8)]{deligne} for a precise statement. 
Our main result on $L$-values is to be viewed from this perspective.

\section{The main result on ratios of critical $L$-values}

\begin{thm}
\label{thm:main}
Let $\sigma_f\in {\rm Coh}(\GL_n, \lambda)$ and $\sigma'_f \in {\rm Coh}(\GL_{n'}, \lambda').$ Assume that $n$ is even and $n'$ is odd. 
Let $m = 1/2 + m_0 \in 1/2 + \Z$ be a half-integer such that both $m$ and $m+1$ are critical for $L(\sigma_f\times \sigma_f^{\prime {\sf v}},\iota, s).$
Assuming the validity of a combinatorial lemma (see below) we have 
$$
\frac{1}{ \Omega(\sigma_f,\iota)^{\epsilon_{m}\epsilon_{\sigma'}}}\frac {\Lambda(\sigma_f\times \sigma_f^{\prime {\sf v}},\iota, m)}
{\Lambda(\sigma_f\times \sigma_f^{\prime {\sf v}},\iota, m+1)} \in
 {\iota(F)}, 
$$
for any $\iota \in \CI(F)$. Moreover, for all $\tau\in \Gal(\bar \Q/\Q)$ 
$$  
\tau\left(\frac{1}{\Omega(\sigma_f,\iota)^{\epsilon_{m}\epsilon_{\sigma'}}}
\frac {\Lambda(\sigma_f\times \sigma_f^{\prime {\sf v}},\iota, m)}
{\Lambda(\sigma_f\times \sigma_f^{\prime {\sf v}},\iota, m+1)}  \right)
\ = \  
\frac{1}{\Omega(\sigma_f,\tau(\iota))^{\epsilon_{m}\epsilon_{\sigma'}}}
\frac {\Lambda(\sigma_f\times \sigma_f^{\prime {\sf v}},\tau(\iota), m)}
{\Lambda(\sigma_f\times \sigma_f^{\prime {\sf v}},\tau(\iota), m+1)}. 
$$
Here 
$\epsilon_{\sigma'}$ is a sign determined by $\sigma'$, $\epsilon_{m} = (-1)^{m_0}$ and
$\Lambda(\sigma_f\times \sigma_f^{\prime {\sf v}},\iota, s)$ is the completed Rankin--Selberg $L$-function.
 \end{thm}

See the main theorem of Harder~\cite{harder-rank-one} for the simplest nontrivial case ($n=2$ and $n'=1$) of the above theorem.

\section{Eisenstein cohomology and sketch of proof of Theorem~\ref{thm:main}}

Consider the group $\tilde G=\GL_N/\Q$ where $N = n+n' \geq 3$ is an odd integer.
Let $P$ (resp., $Q$) be the standard maximal parabolic subgroup of $\tilde G$ 
whose Levi quotient is  $M_P = \GL_n\times \GL_{n^\prime} $ (resp., $M_Q =\GL_{n^\prime}\times \GL_{n }$). We will try to find a 
highest weight  $\tilde \mu,$ such that 
$
 H^{b_n+b_{n'}}_!(S^{M_P}, \E_{\mu, F}) (\sigma_f  \otimes \sigma_f') \ \oplus \ 
 H^{b_n+b_{n'}}_!(S^{M_Q}, \E_{\mu, F}) (\sigma'_f  \otimes \sigma_f)
$
occurs as isotypical summand in the cohomology of the boundary $H^{b_N}(\partial S_{K_f}^{\tilde G},\E_{\tilde \mu}).$ Recall our notation that $b_N = (N^2-1)/4$, hence $b_N = b_n + b_{n'} + \dim(U_P)/2$. Therefore, we need a dominant weight $\tilde \mu$ and a Kostant representative 
$w\in W^P$ (defined as in Borel-Wallach \cite[III.1.2]{borel-wallach}) of length $l(w)=\dim(U_P)/2$  
such that $w\cdot \tilde\mu:= w(\tilde\mu +\tilde\rho)-\tilde \rho= \mu = \lambda+\lambda^\prime.$ 
We believe, having checked it in infinitely many cases ($n=2$ or $n'=1$), that the following assertion is true.

\begin{con}(Combinatorial Lemma)
For a given $\mu=\lambda +\lambda^\prime$, there exists a dominant weight $\tilde\mu$ and a Kostant representative $w \in W^P$ with 
$l(w) = \dim(U_P)/2$ and $w\cdot \tilde\mu = \mu$ if and only if 
$$ 
\frac{(\w+\w')}{2} - p(\mu)+ 1 - \frac{N}{2}
\  \leq \ a(\mu)  \ \leq \ 
 -\frac{(\w+\w')}{2} + p(\mu) -1 - \frac{N}{2}.
$$
(The number of possibilities for $a(\mu)$ is $2p(\mu) - (\w+\w') -1$, which is one less than the total number of critical points.)
\end{con}

Assuming that $\mu$ satisfies the condition in the combinatorial lemma, we know that there is a $\tilde\mu$ such that  
$$
 H^{b_n+b_{n'}}_!(S^{M_P}, \E_{\mu, F}) (\sigma_f  \otimes \sigma_f') \ \oplus \ 
 H^{b_n+b_{n'}}_!(S^{M_Q}, \E_{\mu, F}) (\sigma'_f  \otimes \sigma_f)
 \subset 
 H^{b_N}(\partial S^{\tilde G}, \E_{\tilde \mu}), 
 $$
and it is actually an isotypical subspace. Hence, there is a Hecke-invariant projector $R_{\pi_f}$ to this subspace.  The theory of 
Eisenstein cohomology gives a description of the image of the restriction map 
$$ 
r^* : H^{b_N}(  S^{\tilde G}, \E_{\tilde \mu}) \to H^{b_N}(\partial S^{\tilde G}, \E_{\tilde \mu}).
$$
Our main result on Eisenstein cohomology is the following:

 \begin{thm}
 \label{thm:eisenstein}
The image of $R_{\pi_f}\circ r^*$ is given by 
$$ 
R_{\pi_f}\circ r^*(H^{b_N}(  S^{\tilde G}, \E_{\tilde \mu}) )\otimes_{F,\iota}\otimes \C = 
\left\{\psi  \ + \ 
\frac{C(\mu)}{\Omega(\sigma_f,\iota)^{\epsilon_{\nu_0}\epsilon_{\sigma'}}}
\frac 
{\Lambda^{\rm coh}(\sigma_f\times \sigma_f^{\prime {\sf v}},\iota, \nu_0 )}
{\Lambda^{\rm coh}(\sigma_f\times \sigma_f^{\prime {\sf v}},\iota, \nu_0+1)} 
T^{\rm arith}(\pi_f,\iota)(\psi) \right\},
$$
where $\psi$ is any class in $H^{b_n+b_{n'}}_!(S^{M_P}, \E_{\mu, F}) (\pi_f )$ with $\pi_f = \sigma_f  \otimes \sigma_f'$; 
the operator $T^{\rm arith}(\pi_f,\iota)$ is defined as $T^{\rm arith}(\sigma_f,\iota) \otimes 1_{\sigma_f'}$ after using the K\"unneth-formula; 
$C(\mu)$ is a non-zero rational number; and the point of evaluation is 
$\nu_0 = \frac{\w+\w'}{2}  - a(\mu) - \frac{N}{2}.$ (Note that $\Lambda^{\rm coh}(\sigma_f\times \sigma_f^{\prime {\sf v}},\iota, \nu_0) = 
\Lambda (\sigma_f\times \sigma_f^{\prime {\sf v}},\iota, -N/2).$)
\end{thm}

Theorem~\ref{thm:eisenstein} implies the rationality result stated in Theorem~\ref{thm:main} for $m = -N/2$ because the ratio of $L$-values together with the period is the 
`slope' of a rationally defined map. For an integer $l$, let us change $\sigma$ to $\sigma \otimes \vert \ \vert^l$, then $\lambda$ changes to 
$\lambda - l \cdot {\rm det}$ and 
$a(\mu)$ changes to $a(\mu)-l$, however the possibilities for $l$ are restricted by the inequalities in the Combinatorial Lemma since $\w, \w'$ and $p(\mu)$ do not change. It may be verified using (\ref{eqn:critical}) that as $a(\mu)$ runs through all the possible values it can take as prescribed by 
the Combinatorial Lemma, the pair of numbers $\nu_0$ and $\nu_0+1$ run through all the successive critical arguments; Theorem~\ref{thm:main} follows while using the period relations (\ref{eqn:period-relations}) for $\sigma_f$. The Combinatorial Lemma says that the theory of Eisenstein cohomology allows one to prove a rationality result for a ratio of successive $L$-values exactly when both the $L$-values are critical. 
(See also ~\cite{harder-gln}.)

The condition on $\mu$ imposed by the Combinatorial Lemma has certain strong implications on the situation that underlies Eisenstein cohomology.  
First, using Speh's results (see, for example, \cite[Theorem 10b]{moeglin}) on reducibility for induced representations for $\GL_N(\R)$, one sees that the representation 
${}^{\rm a}{\rm Ind}_{P_{\infty}}^{\GL_N(\R)}(\sigma_{\infty}^{\lambda} \otimes \sigma_{\infty}'^{\lambda'})$ 
of $\GL_N(\R)$ obtained by un-normalized parabolic induction is irreducible. Next, using Shahidi's results \cite{shahidi-duke85} on 
local factors and the fact that $\nu_0$ and $\nu_0+1$ are critical, we deduce that the standard intertwining operator $A_{\infty}$ from the above induced representation to the representation similarly induced from $Q_{\infty}$ is both holomorphic and nonzero at $s = \nu_0$. The choice of bases $\omega_{\pm}$ fixes a basis for the one-dimensional space 
$H^{b_N}(\gl_N, K_\infty^\circ,  {}^{\rm a}{\rm Ind}_{P_{\infty}}^{\GL_N(\R)}(\sigma_{\infty}^{\lambda} \otimes \sigma_{\infty}'^{\lambda'}) \otimes 
E_{\tilde\mu}).$ 
The map induced by $A_{\infty}$ at the level of $(\gl_N, K_\infty^\circ)$-cohomology is then a nonzero scalar. This scalar is a power of $(2 \pi i)$ times 
a rational number $C(\mu)$. The power of $(2 \pi i)$ gives the ratio of $L$-factors at infinity hence giving us a statement for completed $L$-functions, 
and the quantity $C(\mu)$ is expected to be a simple number as was verified for $\GL_3$ by Harder~\cite{harder-rank-one}.


\bigskip

{\small
G\"unter Harder, Max-Planck-Institut f\"ur Mathematik, 
Vivatsgasse 7, 53111 Bonn, Germany. 
E-mail address: {\tt harder@mpim-bonn.mpg.de}

A. Raghuram,
Department of Mathematics,
Oklahoma State University,
401 Mathematical Sciences,
Stillwater, OK 74078, USA. E-mail address: 
{\tt araghur@math.okstate.edu} }

\end{document}